\documentclass[a4paper,12pt]{article}
\usepackage{latexsym}
\usepackage{amssymb}

\makeatletter                                            %This piece serves
\renewcommand{\@begintheorem}[2]{                        %to put numbers
\rm \trivlist \item [\hskip \labelsep {\bf #2\ \ #1.}]   %of theoremlike
                                }                        %environments
\makeatother                                             %in front.

\makeatletter

\newcounter{itemc}[subsection]

\newcommand{\newsubsection}[1]%
{\vspace{\baselineskip}
\noindent\refstepcounter{subsection}\textbf{\thesubsection\ \ #1}}
\newcommand{\ts}{\vspace{\baselineskip}\noindent{\bf Proof.}$\;\;$}
\newcommand{\qed}{{\unskip\nobreak\hfill\hbox{ $\Box$}\par}}
\newcommand{\ZZ}{{\bf Z}}
\newcommand{\QQ}{{\bf Q}}
\newcommand{\RR}{{\bf R}}
\newcommand{\CC}{{\bf C}}

\newcommand{\HH}{{\bf H}}
\newcommand{\PP}{{\bf P}}

\newcommand{\tF}{{\tilde{F}}}
\newcommand{\isog}{{\thickapprox}}

\newcommand{\goth}{{\mathfrak h}}

\begin{document}

\title{Real multiplication on K3 surfaces \\ and Kuga Satake varieties}
\author{Bert van Geemen
\thanks{Universit\`a di Milano, Dipartimento di Matematica,
via Saldini 50, I-20133 Milano, Italia.
{\bf e-mail:} {\it}geemen@mat.unimi.it}}
\date{}
\maketitle
\begin{center}
\abstract{\noindent
The endomorphism algebra of a K3 type Hodge structure is a totally
real field or a CM field. In this paper we give a low brow introduction
to the case of a totally real field.
We give existence results for 
the Hodge structures, for their polarizations and for 
certain K3 surfaces.
We consider the Kuga Satake variety of these Hodge structures
and we discuss some examples.
Finally we indicate various open problems related to the Hodge conjecture.}
\end{center}

\

\

\noindent
A K3 type Hodge structure is a simple, rational, polarized weight 
two Hodge structure $V$ with $\dim V^{2,0}=1$. 
Zarhin \cite{Z} proved that the
endomorphism algebra of a K3 type Hodge structure is either a totally real field or a CM field. 
Conversely, K3 type Hodge structures whose endomorphism algebra
is a given such field exist under fairly obvious conditions.
For the totally real case, see Lemma \ref{no m=2}.

Similar to the case of abelian varieties and their polarized weight one Hodge structures, given a polarization and a totally positive 
endomorphism, one can define a new polarization 
(see Lemma \ref{mult of pols}). 
For a polarized abelian variety,
this follows from the well-known
relation between the Rosati invariant endomorphisms 
and the N\'eron Severi group.
In case the K3 type Hodge structure is a Hodge substructure of the $H^2$
of a smooth surface, it comes with a natural polarization induced by the cupproduct. 
It is then interesting to consider whether the
polarization obtained by means of a totally real element $a$ is also realized as the natural polarization for some other surface $S_a$. 
Thus $H^2(S_a)$ has a Hodge substructure isomorphic 
to the original one, 
but the isomorphism does not preserve the natural polarizations.

It follows easily from general results on K3 surfaces
that, under a condition on the 
dimension of the Hodge structure, 
such K3 surfaces do exist
(see section \ref{twist k3}). 
The isomorphism of Hodge substructures, in combination with the
Hodge conjecture, then leads one to wonder whether there is an algebraic
cycle realizing the isomorphism. 
We discuss some aspects of this question in the last section. 
Mukai \cite{Mukai} proved that Hodge isometries between 
rational Hodge structures of K3 surfaces
are realized by algebraic cycles,
but this deep result does not apply to the general case. 
It does imply that if the endomorphism algebra
of the (transcendental) Hodge structure of the K3 surface
is a CM field, then any endomorphism is induced by an algebraic cycle
on the selfproduct of the surface \cite{pepe}.

We consider the Kuga Satake variety of a K3 type Hodge structure
with real multiplication in sections \ref{kuga satake} and 
\ref{corestriction}. The CM case was already studied in \cite{vG2}.
In particular, we consider the endomorphism algebra of the Kuga Satake variety in the presence of real multiplications on the Hodge structure
and we discuss some examples. We show that the Kuga Satake construction is related to the corestriction of (Clifford) algebras. From this result we obtain a better understanding of previous work of Mumford \cite{M} and
Galluzzi \cite{Ga}.

\section{Hodge structures, polarizations and endomorphisms}

\newsubsection{Hodge structures.}
We recall the basic notions and we refer to \cite{Z}, section 0 and
\cite{vG}, section 1 for further details.

For a $\QQ$-vector space $V$ and a $\QQ$-algebra $R$ we write
$V_R:=V\otimes_\QQ R$. The $R$-linear extension of a $\QQ$-linear
map $f:V\rightarrow W$ is denoted by $f_R:V_R\rightarrow W_R$.

A rational Hodge structure of weight $k$ is a $\QQ$-vector space $V$ with
a decomposition of its complexification:
$$
V_\CC=\oplus_{p+q=k}V^{p,q}\qquad\mbox{such that}\quad
\overline{V^{p,q}}=V^{q,p}
$$
and $p,q\in\ZZ_{\geq 0}$.
Equivalently, a rational Hodge structure of weight $k$ is a
$\QQ$-vector space $V$ with a representation of algebraic groups over $\RR$:
$$
h:U(1)=\{z\in\CC:\;|z|=1\}\longrightarrow GL(V_\RR)
\qquad\mbox{such that}\quad h(z)\sim_\CC diag(\ldots,z^p\bar{z}^q,\ldots)
%\quad(\in GL(V_\CC))
$$
with $p+q=k$, $p,q\geq 0$,
that is the set of eigenvalues of $h(z)$ on $V_\CC$ is a subset of
$\{z^k,\ldots,z^p\bar{z}^q,\ldots, \bar{z}^k\}$.
The subspace of $V_\CC$ on which $h(z)$ acts as $z^p\bar{z}^q$ 
is $V^{p,q}$. Note that $h$ is defined over $\RR$ if and only if
$\overline{V^{p,q}}=V^{q,p}$. We write $(V,h)$, or simply $V$,
for the Hodge structure on $V$ defined by $h$.

\newsubsection{Polarizations.}\label{polar}
(See \cite{Z}, section 0.3.1, 0.3.2 and
\cite{vG}, section 1.7.)
Let $(V,h)$ be a rational Hodge structure of weight $k$. A polarization $\psi$
of $(V,h)$ is first of all a $\QQ$-bilinear map
$$
\psi:V\times V\longrightarrow \QQ\qquad\mbox{such that}\quad
\psi_\RR(h(z)v,h(z)w)=\psi(v,w)
$$
for all $z\in U(1)$, $v,w\in V_\RR$. That is, $\psi_\RR$
must be $U(1)$-invariant.

Let $C:=h(i)\in End(V_\RR)$, $C$ is called the Weil operator.
Then $\psi$ must also satisfy:
$$
\psi_\RR(v,Cw)=\psi_\RR(w,Cv)\qquad(\forall v,w\in V_\RR),
$$
hence $\psi_\RR(\cdot,C\cdot)$ is a symmetric $\RR$-bilinear form.
As $C^2=h(-1)$ and $h(-1)$ has eigenvalues $(-1)^{p+q}=(-1)^k$, 
we have $C^2=(-1)^k$ on $V_\RR$. Using also the symmetry of  $\psi_\RR(\cdot,C\cdot)$, we get:
$$
\psi(v,w)=\psi_\RR(v,w)=\psi_\RR(Cv,Cw)=\psi_\RR(w,C^2v)=
(-1)^k\psi_\RR(w,v)=(-1)^k\psi(w,v),
$$
for  $v,w\in V$.
Thus $\psi$ is symmetric if the weight $k$ of $V$ is even and
alternating if the weight is odd.

Finally one requires that $\psi_\RR(\cdot,C\cdot)$ is positive definite:
$$
\psi_\RR(v,Cv)>0\qquad (\forall v\in V_\RR,v\neq0).
$$
A polarized rational Hodge structure $(V,h,\psi)$
is a rational Hodge structure $(V,h)$ with a polarization $\psi$.

\newsubsection{Endomorphisms.}
A homomorphism $f$ of Hodge structures $(V,h_V)$ and $(W,h_W)$ is a
$\QQ$-linear map whose $\RR$-linear extension intertwines the representations
$h_V$ and $h_W$ of $U(1)$:
$$
f:V\longrightarrow W,\qquad \mbox{such that}\qquad
f_\RR(h_V(z)v)=h_W(z)f_\RR(v)\qquad(\forall z\in U(1),\forall v,w\in V_\RR).
$$
Equivalently, there is an integer $a$ such that the $\CC$-linear extension of $f$ satisfies
$$
f_\CC(V^{p,q})\subset W^{p+a,q+a}
$$
(so we work up to Tate twists, cf. \cite{Z}, 0.3.0, \cite{vG}, 1.6).
In particular kernels and images of homomorphisms of Hodge structures are Hodge substructures, i.e.\ they are rational Hodge structures with
the decomposition induced by the one of $V$ ($W$ respectively).

We write $Hom_{Hod}(V,W)$ for the $\QQ$-vector space of homomorphisms of
Hodge structures and $End_{Hod}(V)=Hom_{Hod}(V,V)$, note that $End_{Hod}(V)$
is a $\QQ$-algebra with product given by the composition.

A Hodge substructure of $(V,h)$ is a subspace $W\subset V$ such that
$h(z)W_\RR\subset W_\RR$ for all $z\in U(1)$. Thus $(W,h_{|W})$ is a rational Hodge structure and the inclusion $W\hookrightarrow V$ is a 
homomorphism of Hodge structures.

A rational Hodge structure $(V,h)$ is said to be simple 
if it does not contain
non-trivial rational Hodge structures. 
If $(V,h)$ is simple then
any non-zero $f\in End_{Hod}(V)$ must be an isomorphism,
i.e.\ $End_{Hod}(V)$ is a division algebra.

\newsubsection{Polarizations on weight two Hodge structures.}
\label{pol weight 2}
Let $(V,h)$ be a rational Hodge structure of weight two.
Define a decomposition of $V_\RR$ (this is actually the eigenspace decomposition for the Weil operator $C$) by:
$$
V_\RR=V_2\oplus V_0,\qquad 
V_2:=V_\RR\cap (V^{2,0}\oplus V^{0,2}),\qquad
V_0:=V_\RR\cap V^{1,1}.
$$
Then $V_2$ is a real vector space which is $h(U(1))$-invariant
and the eigenvalues of $h(z)$ on $V_2\otimes\CC$ are $z^2$ and $\bar{z}^2$. In particular, $C=-1$ on $V_2$.
The subspace $V_0$ is a complementary $h(U(1))$-invariant
subspace of $V_2$ on which $h(U(1))$ acts trivially. In particular
$C=1$ on $V_0$. 

Let $\psi:V\times V\rightarrow \QQ$ be a morphism of Hodge structures.
Then $\psi_\RR$ is $U(1)$-invariant, hence $V_0$ and $V_2$
are perpendicular w.r.t.\ $\psi_\RR$:
$$
\psi_\RR(v,w)=\psi_\RR(Cv,Cw)=-\psi_\RR(v,w)\qquad
(v\in V_0,\,w\in V_2).
$$

Now assume that $\psi$ is a polarization.
Then $\psi_\RR(\cdot,C\cdot)$ is positive definite on $V_\RR$
and  as $C=-1$ on $V_2$, $+1$ on $V_0$,
the symmetric form $\psi_\RR$ is negative definite on $V_2$ 
and positive definite on $V_0$.
In particular, the signature of
$\psi_\RR$ is 
$((d-2e)+,2e-)$ where $d=\dim_\QQ V$, $e=\dim_\CC V^{2,0}$.

\newsubsection{K3 type Hodge structures.}\label{K3 type}
A Hodge structure of K3 type is a
simple, polarized, weight two Hodge structure $(V,h,\psi)$ with
$$
\dim V^{2,0}=1.
$$

\newsubsection{Periods of K3 type Hodge structures.}\label{periods}
Let $(V,\psi,h)$ be a K3 type Hodge structure.
Any non-zero element $\omega\in V^{2,0}$ will be called a period of $(V,h,\psi)$. 
Note that $V^{0,2}=\CC \bar{\omega}$ and $V^{1,1}=(V^{2,0}\oplus V^{0,2})^\perp$, the perpendicular is taken w.r.t.\ the polarization $\psi$. Thus the polarized weight two Hodge
structure $(V,h,\psi)$ is determined by the bilinear form $\psi$,
of signature $((d-2)+,2-)$ where $d=\dim V$ and the period $\omega$.

As $h(z)\omega=z^2\omega$ for all $z\in U(1)$, we get
$$
\psi_\CC(\omega,\omega)=\psi_\CC(h(z)\omega,h(z)\omega)=z^4\psi_\CC(\omega,\omega),\qquad(\forall z\in U(1))
$$ 
hence $\psi_\CC(\omega,\omega)=0$. 
As $\omega+\bar{\omega}\in V_2$ we have 
$0>\psi_\RR(\omega+\bar{\omega},\omega+\bar{\omega})=
2\psi_\CC(\omega,\bar{\omega})$. 
Thus the period satisfies:
$$
\psi_\CC(\omega,\omega)=0,\qquad \psi_\CC(\omega,\bar{\omega})<0,
\qquad (V^{2,0}=\CC \omega).
$$

Conversely, let $V$ be a $\QQ$-vector space of dimension $d$,
with a bilinear form $\psi$ of signature $((d-2)+,2-)$.
Any $\omega\in V_\CC$ which satisfies $\psi_\CC(\omega,\omega)=0$ and $\psi_\CC(\omega,\bar{\omega})<0$
determines a polarized weight two Hodge
structure $(V,h,\psi)$ by
$
V^{2,0}:=\CC \omega
$. 
The following lemma gives a criterion for this Hodge structure to be simple, equivalently, to be of K3 type.

\newsubsection{Lemma.}\label{omega perp}
Let $(V,h,\psi)$ be a weight two Hodge structure with 
$\dim V^{2,0}=1$ and let  $V^{2,0}=\CC\omega$.
The $(V,h,\psi)$ is of
K3 type if and only if  
$\psi_\CC(\omega,v)=0$, with $v\in V$, implies $v=0$.

\ts
We will prove that $V$ has a non-trivial Hodge substructure iff
there is a non-zero $v\in V$ perpendicular to $\omega$.

Let $W\subset V$ be a non-trivial Hodge substructure.
As $W_\RR$ must be invariant under $C=h(i)$, 
it is the direct sum of the perpendicular eigenspaces 
$W_2:=W_\RR\cap V_2$ and $W_0=W_\RR\cap V_0$.

In case $W_2=0$, $W\subset V_0$ and thus 
$\psi_\CC(\omega,w)=0$ for all $w\in W$.
In case $W_2\neq 0$, we get $W_2=V_2$ since $V_2$ is an irreducible (over $\RR$) representation of $U(1)$.
In particular, $\omega,\bar{\omega}\in W_\CC$.  
We consider the subspace
$$
W^\perp:=\{v\in V:\,\psi_\CC(v,w)=0\;\;\forall w\in W\}.
$$
Then $W^\perp$ is also a non-trivial Hodge substructure of $V$
(the $h(U(1))$-invariance of $W^\perp_\RR$ follows from the 
$h(U(1))$-invariance of $W_\RR$ and
$\psi_\RR(h(z)v,h(z)w)=\psi_\RR(v,w)$).
As $\psi_\RR$ is positive definite on $V_2$, it follows that
$W^\perp_\RR\cap V_2=0$.
As before, we find $\psi_\CC(\omega,w)=0$ for all $w\in W^\perp$.

Conversely, given a non-zero $v\in V$ 
with $\psi_\CC(\omega,v)=0$, the $\CC$-linearity of $\psi_\CC$
and the fact that $(\psi_\CC)_{|V_\RR}=\psi_\RR$ 
imply that also $\psi_\CC(\bar{\omega},v)=0$.
Thus $v\in V_2^\perp= V_0$. 
As $h(z)w=w$ for all $z\in U(1)$ and $w\in V_0$,
we conclude that $W=\langle v\rangle$ is a non-trivial Hodge substructure.
\qed

\newsubsection{The transcendental lattice of a K3 surface.} 
\label{trans k3}
Let $S$ be a K3 surface, then $H^2(S,\ZZ)\cong \ZZ^{22}$ and
the natural Hodge structure on $H^2(S,\QQ)$ has
$\dim H^{2,0}(S)=1$ (cf.\ \cite{BPV}). 
The orientation of $S$ gives a natural isomorphism $H^4(S,\ZZ)\cong \ZZ$, so we obtain a cupproduct
$H^2(S,\ZZ)\times H^2(S,\ZZ)\rightarrow \ZZ$. This cupproduct
is an even unimodular bilinear form of signature $(3+,19-)$ on
$H^2(S,\ZZ)$ and (cf.\ \cite{BPV} I.2.7, with $U=H$, VIII.3):
$$
H^2(S,\ZZ)\cong \Lambda_{K3},\qquad 
\Lambda_{K3}=U^3\oplus E_8(-1)^2.
$$

We will assume that $S$ is algebraic, thus $S$ has an ample divisor
with class $h\in H^2(S,\ZZ)\cap H^{1,1}(S)$, in particular 
$h\cup h>0$.
The primitive rational cohomology of $S$ (w.r.t.\ to $h$) is 
$$
h^{\perp}=H^2(S,\QQ)_{prim}=\{v\in H^2(S,\QQ):\,v\cup h =0\}.
$$ 
The inclusion $H^2(S,\QQ)_{prim}\subset H^2(S,\QQ)$ 
defines a rational Hodge 
structure on the primitive cohomology.
The map $\psi_S$, defined by $\psi_S(v,w):=-(v\cup w)$,
gives a polarization on $H^2(S,\QQ)_{prim}$. 

The N\'eron Severi group of $S$ is $NS(S)=H^2(S,\ZZ)\cap H^{1,1}(S)$
(\cite{BPV}, IV, 2.13), note that $NS(S)_\QQ=H^2(S,\QQ)\cap V_0$
is the maximal Hodge structure of type $(1,1)$ contained in $H^2(S,\QQ)$. 
The transcendental lattice of $S$ is defined as
$$
T_S:=NS(S)^\perp\qquad(\subset H^2(S,\ZZ)_{prim}).
$$ 
As $NS(S)\subset V_0$, we get $H^{2,0}(S)\subset T_{S,\CC}$.
The Hodge substructure $T_{S,\QQ}$ of $H^2(S,\QQ)_{prim}$, with the polarization induced by $\psi_S$, is of K3 type. 

The `surjectivity of the period map' (\cite{BPV}, VIII.14)
implies that any $\omega\in \Lambda_{K3}\otimes_\ZZ\CC$
with $\omega\cdot \omega=0$, $\omega\cdot\bar{\omega}>0$
and such that there is a $h\in \Lambda_{K3}$ with $h\cdot h>0$
and $h\cdot \omega=0$, defines an algebraic K3 surface $S$ 
with an isomorphism $H^2(S,\ZZ)\cong \Lambda_{K3}$ whose
$\CC$-linear extension induces an isomorphism 
$H^{2,0}(S)\cong \CC\omega$.

%%%%%%%%%%

\section{Real multiplication for K3 type Hodge structures}

\newsubsection{Endomorphisms of K3 type Hodge structures.}
Zarhin showed that for a Hodge structure
of K3 type $(V,h,\psi)$, the division algebra $End_{Hod}(V)$ is
a (commutative) field which is either totally real,
in which case we write
$End_{Hod}(V)=F$ or it is a CM field $E$, that is $E$ is an imaginary quadratic extension of a totally real field $F$
(\cite{Z}, Theorem 1.5.1).
Moreover, for any polarization $\psi$ of $(V,h)$, one has 
(\cite{Z}, 1.5, Thm 1.5.1)
$$
\psi(av,w)=\psi(v,\bar{a}w),\qquad
(\forall a\in End_{Hod}(V),\;\forall v,\,w\in V)
$$
where $\bar{a}$ is the complex conjugate of $a$, in particular $\bar{a}=a$
for $a\in F$.

\newsubsection{Notation.} \label{notation}
From now on $(V,h,\psi)$ will be a
polarized Hodge structure of K3 type with $F=End_{Hod}(V)$ a
totally real field and
$$
d=\dim_\QQ V,\qquad n=[F:\QQ],\qquad m=n/d=\dim_F V.
$$
%The set of embeddings of $F$ into $\CC$ (whose image will lie in %$\RR$) will be denoted by $S$ and $\epsilon\in S$ will be the 
%embedding defined by $a\cdot \omega=\epsilon(a)\omega$
%where $\omega\in V^{2,0}$ is a period of $(V,h)$.

\newsubsection{Totally real fields.}\label{real fields}
Recall that a finite extension $F$ of $\QQ$ 
is said to be totally real if for any embedding
$\sigma:F\hookrightarrow \CC$ one has $\sigma(F)\subset\RR$.

It is well-known that for any number field $F$ there is an irreducible polynomial $p\in\QQ[X]$ such that $F\cong \QQ[X]/(p)$.
Then $[F:\QQ]=n$, where $n$ is the degree of $p$. 
Let $\alpha_1,\ldots,\alpha_n\in \CC$ be the roots of $p$ in $\CC$.
Then the maps 
$$
\sigma_j:F\cong \QQ[X]/(p)\hookrightarrow \CC,\qquad
\sigma_j:\sum a_iX^i+(p)\mapsto \sum a_i\alpha_j^i,
$$
$j=1,\ldots,n$, 
are the embeddings of $F$ into $\CC$. 
%The set of embeddings of $F$ into $\CC$ will be denoted by 
%$S:=\{\sigma_1,\ldots,\sigma_n\}$.
In particular, $F$ is totally real iff all roots of $p$ are real.

Let $F$ be a totally real field.
An element $a\in F$ 
is called totally positive if $\sigma(a)>0$ for all complex
embeddings $\sigma$ of $F$. For example, $b^2$, for any non-zero
$b\in F$, is totally positive. For any $b\in F$, the element
$k+b\in F$, with $k\in \ZZ_{>0}$, is totally positive
if $k>-\sigma(b)$ for any embedding $\sigma\in S$.

\newsubsection{Splitting over extensions.}\label{split F galois}
To study the action of the field $F=End_{Hod}(V)$ on the 
$\QQ$-vector space $V$ with the bilinear form $\psi$, 
it is convenient to have a decomposition 
into eigenspaces. 
So let $\tF$ be the Galois closure of $F$.
Then $\tF$ is a Galois extension of $\QQ$ which
contains $F$ as a subfield. 
Let 
$$
H:=Gal(\tF/F)\,\hookrightarrow\, G:=Gal(\tF/\QQ),\qquad 
[G:H]=[F:\QQ]=n.
$$
Note that $h(a)=a$ for any $a\in F$ and $h\in H$, thus any coset
$gH$ gives a well-defined embedding $F\hookrightarrow \tF$, $a\mapsto g(a)$.

Let again
$F=\QQ(\alpha)=\QQ[X]/(p)$ with $\alpha=X+(p)\in F$ a root of $p$.
Then $p=\prod_{g\in G/H}(X-g(\alpha))\in \tF[X]$. 
The Chinese remainder theorem gives an isomorphism of rings
$$
F_\tF:=F\otimes_\QQ\tF \,\cong\, 
\prod_{g\in G/H} \tF_g,\qquad
a\otimes t\longmapsto(\ldots,g(a)t,\ldots,)_{g\in G},
$$
where the $F$-algebra
$\tF_g$ is the field $\tF$ on which $F$
acts via the automorphism $g$ of $\tF$: $a\cdot t:=(a\otimes 1)t=g(a)t$
for $t\in \tF_g$.

For $h\in G$, let $\pi_h\in \tF_h\subset F_\tF$ be the idempotent corresponding to the projection on $\tF_h$, so
$\pi_h=(\ldots,(\pi_h)_g,\ldots)\in\prod\tF_g$ with 
$(\pi_h)_h=1$ and $(\pi_h)_g=0$ if $g\neq h$.
Note that $a\cdot \pi_g=g(a)\pi_g$.

Then $V_\tF$ has the decomposition, with $V_g:=\pi_g V$:
$$
V_\tF=\oplus_{g\in G/H}
V_g,\qquad v=\sum_{g\in G/H} v_g 
\qquad\mbox{with}\quad v_g:= \pi_g v.
$$

This is also the decomposition of
$V_\tF$ into eigenspaces for the $F$-action because,
with $a\cdot v=(a\otimes 1)v$ for $a\in F$ and $v\in V_\tF$,
we have
$$
a\cdot v_g=a\cdot \pi_gv_g=g(a)\pi_gv_g=g(a)v_g.
$$
Note that $\dim_\QQ V=\dim_\tF V_\tF=d=nm$ and that
$\dim_\tF V_g=m$ for any $g\in G/H$.

\newsubsection{The Galois action on $V_\tF$.}
The Galois group $G$ acts on 
$V_\tF=V\otimes_\QQ \tF$
via the second factor of the tensor product. 
This action commutes with the one of $F$ on the first factor. 
Under the isomorphism $V_\tF\cong \oplus V_g$, 
this action permutes the eigenspaces $V_g$:
if $a\cdot v=g(a)v$ then for $h\in G$ one has:
$$
a\cdot h(v)=h(a\cdot v)=h(g(a)v)=(hg)(a) h(v)\qquad \mbox{hence}\quad 
h(V_g)= V_{hg}
$$
(use that if $v=\sum_i v_i\otimes t_i$,
then %$a\cdot v=\sum_i (av_i)\otimes t_i$,
$g(a)v=(1\otimes g(a))(\sum_i v_i\otimes t_i)$ etc.).

Let $v\in V=V\otimes 1 \subset V_\tF$ then
$v=h(v)$. Writing $v=\sum v_g$ we get 
$$
v=\sum v_g =\sum h(v_g)\qquad \mbox{hence}\quad h(v_g)=v_{hg},
\qquad(v\in V).
$$
Thus for $v\in V$ we have the decomposition $v=\sum g(v_e)$.
In particular, the composition of the inclusion $V\hookrightarrow V_\tF$
with the projection of $V_\tF\rightarrow V_e$ is an injective $F$-linear map:
$$
V\,\hookrightarrow\, V_e,\qquad v\longmapsto v_e.
$$
This inclusion induces an isomorphism of $\tF$-vector spaces
$V\otimes_F\tF\cong V_e$.

\newsubsection{Lemma.} \label{dec psi}
The $\tF$-bilinear extension  
of the polarization $\psi$ on $V$ will be denoted by $\psi_\tF$,
$$
\psi_\tF:V_\tF\times V_\tF\longrightarrow \tF,\qquad
\mbox{let}\quad \psi_e:V_e\times V_e\longrightarrow \tF
$$
be the restriction of $\psi_{\tF}$ to $V_e\times V_e$.
Let $\Phi$ be the restriction of $\psi_e$ to $V\times V\subset V_e\times V_e$,
$$
\Phi:V\times V\longrightarrow F,\qquad \Phi(v,w):=\psi_e(v_e,w_e).
$$
%be the restriction of $\psi_e$ to $V\times V\subset V_e\times V_e$.

Then $\Phi$ is an $F$-bilinear map and 
for $v=\sum v_g,w=\sum w_g\in V\subset V_\tF$ we have:
$$
\psi(v,w)\,=\,\sum_{g\in G/H} g(\Phi(v_e,w_e))\,
=\,tr_{F/\QQ}(\Phi(v_e,w_e))
$$
where $tr_{F/\QQ}:F\rightarrow \QQ$, $t\mapsto \sum_{g\in G/H} g(t)$ 
is the trace map.

\ts
The polarization $\psi$ on $V$ satisfies $\psi(av,w)=\psi(v,aw)$ 
for all $a\in F$, $v,w\in V$. 
Using the idempotents $\pi_g$ we get:
$$
\psi_\tF(v_g,v_h)=
\psi_\tF(\pi_g v_g,v_h)=
\psi_\tF( v_g,\pi_g v_h),
$$
and $\pi_gv_h=0$ if $g\neq h$.
Thus the eigenspaces $V_g$ are perpendicular w.r.t.\ $\psi_\tF$
and we get
$$
\psi_\tF(v,w)=\sum_{g\in G/H} \psi_g(v_g,w_g),\qquad\mbox{with}\quad 
\psi_g=(\psi_{\tF})_{|V_g\times V_g}:
V_g\times V_g\longrightarrow \tF,
$$
and $v=\sum v_g,w=\sum w_g\in V_\tF$.
Note that $a\psi_g(v_g,w_g)=\psi_g(g(a)v_g,w_g)=\psi_g(v_g,g(a)w_g)$,
thus $\psi_e$ is $\tF$-bilinear.

As $\psi$ is defined over $\QQ$, for $h\in G$ and $v,w\in V_\tF$ we have: 
$$
h(\psi_\tF(v,w))=\psi_\tF(h(v),h(w))
$$
(alternatively, use $v=\sum_i v_i\otimes t_i$ etc.). 
In particular, for $v,w\in V\subset V_e$ and $h\in H$ we get
$\psi_e(v_e,w_e)\in F$.
For
$v=\sum g(v_e),w=\sum g(w_e)\in V\subset V_\tF$ we then have:
$$
%h(\psi_gF(g(v_e),g(w_e))=
%h(\psi_\tF(g(v_e),g(w_e)))=
%\psi_\tF(hg(v_e),hg(w_e))=
%\psi_{hg}(hg(v_e),hg(w_e)),
g(\psi_e(v_e,w_e))=
g(\psi_\tF(v_e,w_e))=
\psi_\tF(g(v_e),g(w_e))=
\psi_{g}(g(v_e),g(w_e)),
$$
hence
$
\psi(v,w)=\psi_\tF(v,w)=\sum \psi_g(g(v_e),g(w_e))=
\sum g(\psi_e(v_e,w_e))$. The lemma now follows from the definition
of $\Phi$.
\qed

\newsubsection{The Hodge group.}\label{zt}
The Hodge group of a K3 type Hodge structure was determined by Zarhin.
For the definition and properties of the Mumford Tate group $MT(V)$ 
and its subgroup, the Hodge group $Hdg(V)=SMT(V)$,
of a Hodge structure $V$ we refer to \cite{Z}, 0.3.1 and \cite{BG}.
Both are algebraic subgroups, defined over $\QQ$,
of $GL(V)$ and $SL(V)$ respectively.

%The inclusion $V=V\otimes 1\hookrightarrow V_\tF= \oplus V_g$
%embeds the $\QQ$-vector space $V$ diagonally into $\oplus V_g$. 
%The inclusion followed by the projection on $V_e$
%is an isomorphism $V\rightarrow V_e$ of $F$-vector spaces. Let 
%$$
%\Phi:V\times V\longrightarrow F
%$$
%be the $F$-bilinear map which corresponds to $\psi_e$ under this %isomorphism.

Let $SO(V,\Phi)$ be the special orthogonal group of the bilinear form $\Phi$ on the $F$-vector space $V$ (defined in Lemma \ref{dec psi}), viewed as an algebraic group over $\QQ$.
Then $SO(V,\Phi)\cong SO(V_e,\psi_e)$.
For a $\QQ$-algebra $R$, the group of $R$-valued points of $SO(V,\Phi)$ is:
$$
SO(V,\Phi)(R)=\{A\in SL(V_R):\,aA=Aa,\quad\Phi_R(Av,Aw)=\Phi_R(v,w),\;\;
\forall a\in F,\;\forall v,w\in V_R\,\}.
$$ 

\newsubsection{Theorem}(\cite{Z}, Thm. 2.2.1) 
\label{real Hodge}
Let $(V,h,\psi)$ be a K3 type Hodge structure with $F=End_{Hod}(V)$
a totally real field, then
$$
Hdg(V)=SO(V,\Phi),\qquad SO(V,\Phi)(\RR)\cong SO(2,m-2)\times SO(m,\RR)^{n-1},
\qquad 
$$
and $
SO(V,\Phi)(\CC)\cong SO(m,\CC)^n.
$
The representations of these Lie groups on the $d=nm$-dimensional
vector spaces $V_\RR$, $V_\CC$ respectively, are the direct
sum of the standard representations of the factors.

\ts
With our definition of $\Phi$, Lemma \ref{dec psi} shows that
$\psi(v,w)=tr_{F/\QQ}(\Phi(v,w))$ for all $v,w\in V$,
and this is also Zarhin's definition of $\Phi$ (\cite{Z}, 2.1, with $e=1$).
Thus $Hdg(V)=SO(V,\Phi)$ by \cite{Z}, Thm 2.2.1.

Similar to the various decompositions in \ref{split F galois}, one
has $F_\RR\cong \oplus_{\sigma\in S}\RR_\sigma$ and
$$
V_\RR=\oplus_{\sigma\in S}V_{\sigma},\qquad v=(\ldots,v_\sigma,\ldots )\qquad\mbox{and}\quad
a\cdot v=(\ldots,\sigma(a)v_\sigma,\ldots)_{\sigma\in S}
$$
for $a\in F$. Thus this is also the decomposition of
$V_\RR$ into eigenspaces for the $F$-action.

To obtain these from the decomposition of $F_\tF$,
choose an embedding $\tilde{\epsilon}:\tF\hookrightarrow \RR$, 
for example one which extends
$\epsilon:F\hookrightarrow \RR$, where 
$a\cdot \omega=\epsilon(a)\omega$ and $V^{2,0}=\CC\omega$.
(We are only interested in the action of $F$ and $\RR$ on $F_\RR$, so
the choice of the extension of $\epsilon$ does not matter).
Then $\RR$ becomes an $\tF$-module which we
denote by $\RR_\epsilon$
%$a\otimes u=\epsilon(a)u$ for $a\in F,u\in \RR$
and we get
$$
F_\RR:=F\otimes_\QQ\RR\cong(F\otimes_\QQ\tF)\otimes_\tF\RR_\epsilon
\cong\, (\oplus_{g\in G/H} \tF_g)\otimes_\tF\RR_\epsilon\,\cong
\oplus_{\sigma\in S}\RR_\sigma
$$
where the bijection between the set of complex embeddings 
$S$ of $F$ and $G/H$ is given by 
$\sigma(a)=\tilde{\epsilon}(g(a))$, for all $a\in F$.
The idempotent in $F_\RR$ which corresponds to the projection on $\RR_\sigma$ will be denoted by $\pi_\sigma$ and $V_\sigma=\pi_\sigma V$. 

The fact that $\Phi$ is $F$-bilinear implies that its
$\RR$-bilinear extension $\Phi_\RR$, 
with values in $F_\RR\cong\prod\RR_{\sigma}$, is $F_\RR$-bilinear. 
Thus $\Phi_\RR=(\ldots,\Phi_\sigma,\ldots)_{\sigma\in S}$ with
$\Phi_\sigma:V_\sigma\times V_\sigma\rightarrow \RR_\sigma$.

Any $F$-linear endomorphism $A$ of $V$ extends $\RR$-linearly to
an $F_\RR$-linear endomorphism $A_\RR$ of $V_\RR$. In particular,
$A_\RR$ commutes with the idempotents $\pi_\sigma\in F_\RR$, 
hence for $v_\sigma\in V_\sigma$ we get $A_\RR v_\sigma=A_\RR\pi_\sigma v_\sigma=\pi_\sigma A_\RR v_\sigma$
hence $A_\RR v_\sigma\in\pi_\sigma V=V_\sigma$. 

Since the elements of $SO(V,\Phi)$ are $F$-linear, we get:
$
SO(V,\Phi)_\RR \subset \prod_{\sigma\in S}GL(V_\sigma).
$
Here $SO(V,\Phi)_\RR=SO(V,\Phi)\times_\QQ\RR$ is the
algebraic group over $\RR$ obtained by extension of scalars
(that is, by base change) from the algebraic group $SO(V,\Phi)$ over $\QQ$.
As $(V,\Phi)\cong (V_e,\psi_e)$ and the Galois group $G$ permutes the
$V_g$, we obtain the following isomorphism of algebraic groups:
$$
SO(V,\Phi)_\RR\cong \prod_{\sigma\in S}SO(V_\sigma,\Phi_\sigma),
$$
and the representation is factorwise on $V_\RR=\oplus V_\sigma$.

By definition of $\epsilon$ we have $V^{2,0}\subset V_{\epsilon,\CC}:=V_\epsilon\otimes_\RR\CC$.
Since $h$ commutes with $F$, each $V_\sigma$ is a real Hodge
structure, hence also $V^{0,2}\subset V_{\epsilon,\CC}$:
$$
V^{2,0}\oplus V^{0,2}\subset V_{\epsilon,\CC}:=V_\epsilon\otimes_\RR\CC
\qquad\mbox{hence}\quad V_2\subset V_\epsilon.
$$

Recall from \ref{pol weight 2} that $V_\RR=V_2\oplus V_0$ and that
$\psi_\RR$ is negative definite on 
$V_2$ and positive definite on $V_0$. 
As $V_\sigma\subset V_0$ if $\sigma\neq \epsilon$,
we get $SO(V_\sigma,\Phi_\sigma)_\RR\cong SO(m)_\RR$ 
if $\sigma\neq \epsilon$ and
$SO(V_\epsilon,\Phi_\epsilon)_\RR\cong SO(m-2,2)_\RR$.
Extending scalars to $\CC$ and taking $\CC$-points, we get $SO(V,\Phi)(\CC)\cong SO(m,\CC)^n$.
\qed

%%%%%%%%%%%%%%

\section{The existence of Hodge structures of K3 type with real multiplication.}

\newsubsection{} We prove an easy existence result for 
K3 type Hodge structures with real multiplication.
Next we apply results of Nikulin on embeddings of lattices
and the surjectivity of the period map 
to show the existence of K3 surfaces whose transcendental lattice 
has real multiplication. 

It would be interesting to have geometrical 
(and not just `Hodge theoretic') examples of such surfaces. 
A first step in this direction is done in \ref{double cover}, where real multiplication for
a certain `geometric' four dimensional family of K3 surfaces 
is studied.

\newsubsection{Lemma.}\label{no m=2}
Let $F$ be a totally real field with $[F:\QQ]=n$.
Then for any $m\in \ZZ_{\geq 3}$ there exist K3 type Hodge
structures $(V,h,\psi)$ with $End_{Hod}(V)=F$ and $\dim_FV=m$.

However, there are no such K3 type Hodge structures 
with $\dim_F V\leq 2$.

\ts
We use the notation from the proof of Theorem \ref{real Hodge}.
In case $m=1$ we would have $1=\dim V_\epsilon$, which contradicts
the fact that the two dimensional subspace
$V_2$ is a subspace of $V_\epsilon$. 

If $m=2$ and
$End_{Hod}(V)$ were equal to $F$, then by Zarhin's theorem
we would have $Hdg(V)(\CC)\cong SO(2,\CC)^n$. 
A basic property of $Hdg(V)$ is that $End_{Hdg}(V)=End_{Hod}(V)$.
Note that $SO(2,\CC)\cong \CC^*$ and its standard representation on $\CC^2$ is equivalent to $t\mapsto diag(t,t^{-1})$. 
With this action,
$End_{\CC^*}(\CC^2)$, the endomorphisms commuting with $\CC^*$,  
consists of the diagonal matrices in $End(\CC^2)$. Hence 
$$
End_{Hod}(V)_\CC\cong (End_{\CC^*}(\CC^2))^n
\cong \CC^{2n},
$$
and thus $\dim_\QQ End_{Hod}(V)=2n$
which contradicts that $End_{Hod}(V)=F$.

Fix an integer $m\geq 3$.
It remains to show that there exist K3 type Hodge structures
with $\dim _FV=m$.
Let $V=F^m$ and choose an embedding $\epsilon:F\hookrightarrow\RR$.
Using the isomorphism $F_\RR\cong \prod\RR_\sigma$
it is easy to see that there are $a_i\in F$ such that:
$$
\epsilon(a_1)<0,\quad\epsilon(a_2)<0,\quad \epsilon(a_j)>0
\qquad \mbox{and}
\qquad
\sigma(a_i)>0\quad\mbox{if}\quad \sigma\neq\epsilon,
$$
with $3\leq j\leq n$, $1\leq i\leq n$.
We define an $F$-bilinear form:
$$
\Phi:V\times V\longrightarrow F,
\qquad
\Phi(x,y)=\sum_{k=1}^m a_kx_ky_k.
$$
Then $\Phi$ induces the bilinear form defined by 
$\sum_{k=1}^m \sigma(a_k)x_ky_k$ on $V_\sigma\cong \RR^m$.
Note that $\Phi_\epsilon$ has signature $((m-2)+,2-)$
and $\Phi_\sigma$ is positive definite if $\sigma\neq \epsilon$.
Thus the signature of the $\QQ$-bilinear form $\psi:=tr(\Phi)$ 
on $V$ is $((nm-2)+,2-)$.

To define a K3 type Hodge structure, with polarization $\psi$, on
$V$ it suffices to give a period $\omega\in V_\CC$
(cf.\ \ref{periods}) such that $\psi_\CC(\omega,v)\neq 0$ for all 
non-zero $v\in V$ (Lemma \ref{omega perp}). 
Since %we want $F\subset End_{Hod}(V)$ 
$\psi_\RR<0$ on $V_2$, we must 
choose $\omega\in V_{\epsilon,\CC}:=V_\epsilon\otimes_\RR\CC$.
As $\omega$ and $\lambda\omega$, $\lambda\in \CC-\{0\}$ define
the same Hodge structure, we consider
$$
{\mathcal D}:=\{[\omega]\in \PP( V_{\epsilon,\CC}):\;
\psi_\CC(\omega,\omega)=0,\quad \psi_\CC(\omega,\bar{\omega})<0\,\},
$$
which is a non-empty open subset in a quadric in a complex projective space of dimension $m-1\geq 2$.
Any non-zero $v\in V$ defines a hyperplane 
$$
H_v:=\{[w]\in \PP( V_{\epsilon,\CC}):\;
\psi_\CC(v,w)=0\,\}\qquad
\subset \PP( V_{\epsilon,\CC}).
$$
As an open subset of a quadric is not contained in a hyperplane,
$H_v\cap {\mathcal D}$ is an analytic subset of codimension
at least one in ${\mathcal D}$. As $V$ is a countable set, 
we get ${\mathcal D}\neq \cup_v (H_v\cap {\mathcal D})$, where the union
is over the non-zero $v\in V$. Hence there is an 
$\omega\in{\mathcal D}$ which defines a simple Hodge structure with
$V^{2,0}=\CC\omega$. This construction shows that such
(integral) Hodge structures have $m-2$ moduli.

As $\omega\in V_{\epsilon,\CC}$ and $a\in F$ acts by scalar multiplication by $\epsilon(a)\in \RR$ on $V_{\epsilon,\CC}$, it follows that $aV^{2,0}\subset V^{2,0}$ and, taking the complex conjugates, $aV^{0,2}\subset V^{0,2}$. As $V^{1,1}=(V^{2,0}\oplus V^{0,2})^\perp$ and $\psi_\CC(av,w)=\psi_\CC(v,aw)$ we get $aH^{1,1}\subset H^{1,1}$, so $F\subset End_{Hod}(V)$.

In case $F\neq End_{Hod}(V)$, the K3 type Hodge structure $V$ has
$End_{Hod}(V)=F'$, with $F'$ an extension of the field $F$.
Multiplication by $b\in F'$ defines an $F$-linear map 
$V\rightarrow V$, and thus each eigenspace $V_\sigma$ for the $F$-action is mapped into itself.
The splitting of $V_\CC$ into $F'$ eigenspaces thus splits each
$$
V_{\sigma,\CC}=\,\oplus_\rho\, V_{\rho}
$$
where the $\rho:F'\rightarrow \CC$ are the embeddings of $F'$ which restrict to $\sigma$ on $F\subset F'$.
In particular, each $V_\rho$ has dimension at most 
$(\dim_\CC V_{\sigma,\CC})/2$.
As $F'\in End_{Hod}(V)$ we must have $\omega\in V_{\rho}$ for some
$\rho$ which extends $\epsilon$. 
Since the $V_\rho$ are eigenspaces of elements $b\in F'\subset End_F(V)$ and the set $End_F(V)$ is countable, 
we conclude that the general $\omega\in {\mathcal D}$
defines a K3 type Hodge structure $V$ with $End_{Hod}(V)=F$.
\qed

\newsubsection{Proposition.}\label{exist rm}
Given a totally real number field $F$ and an integer $m\geq 3$ such that
$m[F:\QQ]\leq 10$, there exist $(m-2)$-dimensional families
of K3 surfaces such that $F= End(T_S)$ for the general surface
$S$ in the family.

\ts
Let $(V,h,\psi)$ be a K3 type Hodge structure with 
$End_{Hod}(V)=F$ and period $\omega$. 
Choose a free 
$\ZZ$-module $T\subset V$ of rank $d=\dim_\QQ V$ 
such that $\psi$
is integer valued on $T\times T$.
Theorem 1.10.1 of \cite{Ni}
shows that there is a primitive embedding of lattices 
$T\hookrightarrow \Lambda_{K3}$. 
The surjectivity of the period map implies that 
$\omega\in T\otimes_\ZZ\CC=V_\CC\subset \Lambda_{K3,\CC}$ 
defines a K3 surface $S$ with 
$T_S\cong T$ as integral polarized Hodge structures.
The proof of Lemma \ref{no m=2} shows that there are $m-2$ moduli.
\qed

\newsubsection{Example.}\label{double cover}
The minimal model of a double cover of $\PP^2$ branched over six lines is
a K3 surface. The general surface $S$ in this four dimensional family
has (cf.\ \cite{matsumoto} 0.3):
$$
T_S\cong U(2)^2\oplus <-2>^2,\qquad {\rm hence}\quad
T_{S,\QQ}\cong\, <1>^2\oplus <-1>^4
$$
where $U(2)$ is the lattice $\ZZ^2$ with quadratic form $4x_1x_2$, which is isomorphic over $\QQ$ to
$y_1^2-y_2^2$ (put $x_1=y_1+y_2$, $x_2=(y_1-y_2)/2$), 
%which is denoted by $<1>\oplus <-1>$,
and $<-2>^2$ is $\ZZ^2$ with quadratic form $-2x_1^2-2x_2^2$ 
which is isomorphic to $<-1>^2$ (put $x_1=(y_1+y_2)/2$, $x_2=(y_1-y_2)/2$).

We will show that there are one parameter families
of K3 surfaces such that $T_S\cong T$ and $End_{Hod}(T_S)$ 
is a real quadratic field for the general $S$ in the family.

We consider the vector space $V=\QQ^6=(\QQ^2)^3$ with the following
bilinear form $\psi$:
$$
(V,\psi)\,=\,
(\QQ^2,Q_{1})\oplus (\QQ^2,Q_{2})\oplus (\QQ^2,Q_{3}),\qquad
Q_i:=\left(\begin{array}{rr}1&0\\0&r_i\end{array}\right),
$$
so we identify the bilinear form with the symmetric matrix which defines it and $r_i\in\QQ$ are to be chosen later.
Next we want to consider the $a\in End(V)$ 
such that $\psi(ax,y)=\psi(x,ay)$ for all $x,y\in V$.
We will restrict ourselves to those which  
preserve the direct sum decomposition. Thus
$a=block(A_1,A_2,A_3)$ with $A_i\in End(\QQ^2)$ and
${}^tA_iQ_i=Q_iA_i$, so we have to consider the matrix equation
$$
\left(\begin{array}{rr}a&c\\b&d\end{array}\right)
\left(\begin{array}{rr}1&0\\0&r\end{array}\right)=
\left(\begin{array}{rr}1&0\\0&r\end{array}\right)
\left(\begin{array}{rr}a&b\\c&d\end{array}\right),
\qquad\mbox{let}\quad
A_{e,c}:=\left(\begin{array}{rr}e&cr\\c&-e\end{array}\right),
$$
then the solutions to the matrix equation are $A=A_{e,c}+\lambda I$,
note that $A_{e,c}$ is just the `traceless part' of $A$.
As $A^2_{e,c}=(e^2+rc^2)I$, $\QQ(A_{e,c})\cong \QQ(\sqrt{e^2+rc^2})$ is a quadratic extension $\QQ$ if
$e^2+rc^2$ is not a square in $\QQ$.

Now we consider the case $r_1=r_2=-1$ and $r_3=+1$, 
so $(V,\psi)\cong T_{S,\QQ}$ for a general $S$ as above.

Assume that $d\in\ZZ_{>0}$ is odd, square free and that $d=e^2+c^2$ 
for some integers $c,e$. Write $d=2d'+1$, then $d=(d'+1)^2-(d')^2$. 
Hence if we define
$$
a=(A_1,A_2,A_3),\qquad\mbox{with}\quad 
A_1=A_2=A_{d'+1,d'},\quad A_3=A_{e,c},\quad(c^2+e^2=d=2d'+1)
$$ 
then $a^2=d$ so we have an action of the 
real field $F=\QQ(a)\cong \QQ(\sqrt{d})$ on $V$ and the elements
of $F$ are self adjoint for the bilinear form $\psi$.

It is easy to see that one eigenspace for the $F$-action on $V_\RR$
is positive definite for $\psi$ and the other, call it $V_\sigma$, has signature $(1+,2-)$. 
Next we choose $T\cong\ZZ^6\hookrightarrow \QQ^6$ 
such that $(T,\psi)\cong T_S=U(2)^2\oplus <-2>^2$
so the lattice $(T,\psi_T)$ is isometric
to the lattice of double cover of $\PP^2$ branched over 6 lines.
Choosing a general $\omega\in V_{\sigma,\CC}$ with 
$\psi_{T,\CC}(\omega,\omega)=0,
\psi_{T,\CC}(\omega,\bar{\omega})>0$
defines a polarized integral Hodge structure on $(T,\psi_T)$
with $End_{Hod}(T_{\QQ})=F$ (cf.\ the proof of Lemma \ref{no m=2}). 
Thus  we obtain a one parameter family of Hodge
structures on $(T,\psi_T)$ with $End_{Hod}(T_\QQ)=F$ for the general
member of the family.

\section{Twisting the polarization}

\newsubsection{Real multiplication and polarizations.}
Let $(V,h,\psi)$ be a K3 type Hodge structure.
Let $B_1(V)\subset V^*\otimes V^*$ be the subspace of 
$\QQ$-bilinear maps $\phi:V\times V\rightarrow \QQ$ such that $\phi_\RR$ is $h(U(1))$-invariant.

The isomorphism
$End(V)= V^*\otimes V\cong V^*\otimes V^*$, given by the
isomorphism $V\rightarrow V^*$ defined by $\psi$,
defines an isomorphism
$$
End_{Hod}(V)\longrightarrow B_1,\qquad a\longmapsto \psi_a,
\quad \mbox{with}\quad \psi_a(v,w)=\psi(av,w).
$$
The bilinear form $\psi_a$ is symmetric iff $\QQ(a)$ is a totally
real field (use $\psi(av,w)=\psi(v,\bar{a}w)$).

We now consider when the bilinear form $\psi_a$ 
gives a polarization on $V$.

\newsubsection{Lemma.}\label{mult of pols}
Let $(V,h,\psi)$ be a Hodge structure of K3 type such that
$F=End_{Hod}(V)$ is a totally real field. 
Then $\psi_a$ is a polarization of the Hodge structure 
$(V,h)$ if and only if
$a\in F$ is totally positive.

\ts
As $\psi_\RR$ is $U(1)$-invariant and
$a$ commutes with $h(U(1))$
it follows that $\psi_{a,\RR}$ is also $U(1)$-invariant.
As $C=h(i)$ we get $aC=Ca$ and as $\psi(av,w)=\psi(v,aw)$
it follows that $\psi_{a,\RR}(\cdot,C\cdot)$ is symmetric.
Thus $\psi_a$ is a polarization iff 
$\psi_{a,\RR}(v,Cv)>0$ for all non-zero $v\in V_\RR$.

Using the decomposition $V_\RR=\oplus V_\sigma$ (cf.\ the proof of \ref{real Hodge}
and \ref{dec psi}), 
for $v=(\ldots,v_\sigma,\ldots)\in V_\RR$ we have:
$$
\psi_{a,\RR}(v,Cv)=
\sum_{\sigma\in S} \psi_\sigma(a\cdot v_\sigma,Cv_\sigma)=
\sum_{\sigma\in S} \sigma(a)\psi_\sigma(v_\sigma,Cv_\sigma).
$$
As $\psi$ is a polarization,  $\psi_\sigma(\cdot,C\cdot)$ is positive definite
for any $\sigma\in S$.
Thus $\psi_a$ is a polarization iff 
$\sigma(a)>0$ for all $\sigma\in S$ iff 
$a$ is totally positive.
\qed 

\newsubsection{Example.}
Let $(V,h,\psi)$ and $F$ be as in the Lemma (or equivalently, as in \ref{notation}).
In case $a=b^2$ for some $b\in F$ we have $\psi(av,w)=\psi(bv,bw)$.
Thus the map $B:V\rightarrow V$ given by multiplication by $b$,
$Bv:=bv$ is an isometry between $V,\psi_a)$ and $(V,\psi)$.
As $b$ is a map of Hodge structures, $B$ is a Hodge isometry from
$(V,h,\psi_b)$ to $(V,h,\psi)$.

To find examples where there is no Hodge isometry between the K3 type Hodge structure and its twist, the determinant of a bilinear form is useful.

\newsubsection{Determinants and discriminants.}\label{det and dis}
Let $\psi,\psi_a$ be two bilinear forms on a $\QQ$ vector space $V$.
Choose a basis of $V$.
If $Q,Q_a$ are the symmetric matrices defining the 
symmetric bilinear forms $\psi,\psi_a$ respectively, and
$B$ is (the matrix of) an isometry between $(V,\psi_a)$ and $(V,\psi)$,
we must have ${}^tBQB=Q_a$, in particular, $\det(Q_a)=\det(B)^2\det(Q)$. We will write $\det(\psi)$ for (the class
of) $\det(Q)$ in $\QQ^*$ (modulo the subgroup of squares).
This gives a well-known invariant (often called disciriminant) of a quadratic space. 
%We will only
%need that if $\det(Q)/\det(Q_a)$ is not a square in $\QQ$ then 
%there cannot exist an isometry between 
%$(V,\psi)$ and $(V,\psi_a)$.

If $F$ is a finite extension of $\QQ$, the discriminant of $F$
is the rational number $D_F$, well defined up to squares in $\QQ$, 
defined as
$
D_F:=\det(tr(e_ie_j))
$
where the $e_i\in F$ are a $\QQ$-basis of $F$. 

\newsubsection{Lemma.}
\begin{enumerate}
\item{} Let $\psi=tr(\Phi)$ as in Lemma \ref{dec psi} and $m=\dim_FV$. Then
$$
\det(\psi)\,=\,D_F^mN(\det(\Phi)),
$$
where $m=\dim_FV$.

\item{} For $a\in F$, let $\psi_a(\cdot,\cdot):=\psi(a\cdot,\cdot)$. Then
$$
\det(\psi_a)=N(a)^m\det(\psi),
$$
where $N(a):=\prod_{\sigma\in S}\sigma(a)$ is the norm of $a$. 
\end{enumerate}
\ts
Choose an $F$-basis of $V$ for which $\Phi$ is diagonal:
$\Phi(x,y)=\sum a_k x_k y_k$. Then
$$
\det(\psi)=\prod \det(\psi^{(k)})\qquad\mbox{ with}\qquad
\psi^{(k)}:F\times F\longrightarrow \QQ,\quad \psi^{(k)}(x,y)=tr(a_kxy).
$$
Let $e_1,\ldots,e_n$ be a $\QQ$-basis of $F$ and let 
$S=\{\sigma_1,\ldots,\sigma_n\}$.
As $D_F=\det(\sigma_i(e_j))^2$ (cf.\ \cite{samuel}, Proposition II.3),  
one finds
$$
\begin{array}{rcl}
\det(tr(ae_ie_j))&=&\det(\sum_k\sigma_k(a)\sigma_k(e_i)\sigma_k(e_j))\\
&=&\det\left(\sigma_k(a)\sigma_k(e_i)\right)
\det\left(\sigma_k(e_j)\right)\\
&=&(\prod_k\sigma_k(a))\det\left(\sigma_k(e_i)\right)
                       \det\left(\sigma_k(e_j)\right)\\
&=&N(a)D_F.
\end{array}
$$
Thus $\det(\psi)=D_F^m\prod_kN(a_k)=D_F^mN(\det(\Phi))$.
%To compute $\det(\psi_a)$ in terms of $a$ and $\det(\psi)$ 
To compute $\det(\psi(a)$ use that
$$
\psi_a(v,w)=\psi(av,w)=\psi_\RR(av,w)=
\oplus_{\sigma\in S} \sigma(a)\Phi_\sigma(v_\sigma,w_\sigma),
$$
thus, with $m=\dim_FV=\dim_\RR V_\sigma$, we get:
$$
\det(\psi_a)=\prod_{\sigma\in S} \sigma(a)^m\det(\Phi_\sigma)=
N(a)^m\det(\psi).
$$
\qed

\newsubsection{Examples.}
In case $F\cong \QQ(\sqrt{d})$, with $d$ square free, and $m=\dim_FV$ is odd, 
it is easy to produce
examples of totally positive $a\in F$ such that 
$(V,\psi)$ and $(V,\psi_a)$ are not isometric. 
In fact, $d\pm \sqrt{d}>0$,
so $a=d+\sqrt{d}$ is totally positive. As $N(a)=d^2-d=d(d-1)$ and
$d$ is square free, $N(a)$ is not a square in $\ZZ$, hence $\det(\psi_a)/\det(\psi)$ is not a square in $\QQ$.

%It follows from classical results
%in number theory that there exist totally positive $a\in F$ for which %$N(a)$ is not a square. Thus if $m$ is odd, there exist
%polarized Hodge structures of K3 type $(V,h,\psi)$ and $(V,h,\psi_a)$ %which are not Hodge isometric.

\newsubsection{Twisting K3 surfaces with real multiplication.}
\label{twist k3}
Let $(V,h,\psi)$ be a Hodge structure of K3 type, with $\dim V\leq 11$.
Then for any totally positive $a\in F=End_{Hod}(V)$ we obtain the polarized 
Hodge structure of K3 type $(V,h,\psi_a)$.

Results of Nikulin and the surjectivity of the period map
(cf.\ the proof of Proposition \ref{exist rm}) imply that
there exist K3 surfaces $S,S_{a}$ such that 
$$
(V,h,\psi)\cong T_{S,\QQ},\qquad 
(V,h,\psi_a)\cong T_{S_{a},\QQ},
$$
where the polarizations on the right hand sides are induced by 
(minus) the cup product in the corresponding surface.
We will call $S_{a}$ a real twist of $S$. Note that there
are isomorphisms of rational Hodge structures 
$(V,h)\cong T_{S,\QQ}\cong T_{S_{a},\QQ}$, but in
general there is no isomorphism of
polarized Hodge structures between $T_{S,\QQ}$ and $T_{S_{a},\QQ}$.

%%%%%%%%%%%%%%%%%%%

\section{The Kuga-Satake variety}\label{kuga satake}

\newsubsection{The Kuga Satake construction.}
We briefly recall the construction of `the' Kuga Satake variety,
actually an isogeny class of abelian varieties,
of a rational, polarized,
Hodge structure $(V,h,\psi)$ of weight two with
$\dim V^{2,0}=1$ (cf.\ \cite{vG}, section 5).

%Let $Q:V\rightarrow \QQ$ be the non-degenerate quadratic form
%defined by $Q(v)=\psi(v,v)$.
The Clifford algebra $C(\psi)$ 
%of the quadratic space $(V,\psi)$
is the quotient of the tensor algebra 
$T(V)=\oplus_n V^{\otimes n}$ by the two sided ideal generated
by $v\otimes v-\psi(v,v)$, where $v$ runs over $V$.
The dimension of $C(\psi)$ is $2^d$, where $d=\dim V$.
The Clifford algebra has a subalgebra $C^+(\psi)$ of dimension $2^{d-1}$, 
the quotient of $\oplus_m V^{\otimes 2m}$,
called the even Clifford
algebra.

The Hodge decomposition of $V$ defines 
the subspace $V_2\subset V_\RR$. Choosing a basis $f_1,f_2$ of $V_2$
such that $\langle f_1+if_2\rangle=V^{2,0}$  and $\psi_\RR(f_1,f_1)=-1$ 
defines an element $f_1f_2\in C^+(\psi)_\RR$. 
Multiplication by $f_1f_2$ defines a map $J:c\mapsto f_1f_2c$, 
which is a complex structure on $C^+(\psi)_\RR$, $J^2=-I$ 
(cf.\ \cite{vG}, Lemma 5.5). This defines a weight one Hodge structure
$h_s$ on $C^+(\psi)$ by
$$
h_s:U(1)\longrightarrow GL(C^+(\psi)_\RR),\qquad a+bi\longmapsto a+bJ
$$
for $a,b\in\RR$.
The choice of
$e_1,e_2\in V$ with $\psi(e_1,e_1)>0,\psi(e_2,e_2)>0$ and $\psi(e_1,e_2)=0$
determines a Riemann form, i.e.\ an alternating form 
$$
E:C^+(\psi)\times C^+(\psi)\longrightarrow \QQ\qquad
\mbox{such that}\quad
E(Jx,Jy)=E(x,y),\quad E(x,Jx)>0
$$
for all $x,y\in C^+(\psi)_\RR$ (cf.\ \cite{vG}, Proposition 5.9).
 The complex structure $J$ is
uniquely determined by $(V,h,\psi)$, but
the polarization, as constructed in \cite{vG}, 5.7, depends on the
choice of a positive 2-plane in $V$ and is not unique in general.

Thus $(V,h,\psi)$ defines a polarized rational Hodge structure of weight one $(C^+(\psi),h_s,E)$. 
Each abelian variety $A$ in the isogeny class of abelian varieties
of dimension $2^{d-2}$ defined by 
$(C^+(\psi),h_s,E)$ will be called a Kuga Satake variety for 
$(V,h,\psi)$, so $A$ is characterized by an isomorphism of Hodge structures
$H^1(A,\QQ)\cong (C^+(\psi),h_s)$.

\newsubsection{The endomorphism algebra of the Kuga Satake variety.}
Let $(V,h,\psi)$ be a general polarized Hodge structure of K3 type, more precisely, assume that $MT(V)=GO(\psi)$
where $GO(\psi)=\{g\in GL(V):\psi(gv,gw)=\lambda_g \psi(v,w)\,\}$. 
Then 
$$
MT(A)=CSpin(\psi),\qquad End(A)_\QQ\cong C^+(\psi)
$$ 
(\cite{vG}, Proposition 6.3.1, Lemma 6.5), where 
$CSpin(\psi)$ is the Spin group of $(V,\psi)$. 
There is an isomorphism of complex Lie algebras 
$Hdg(A)_\CC\cong so(d)_\CC$,
$d=\dim V$, and $H^1(A,\CC)$ is a direct sum of copies of the Spin representation $S(d)$ of $so(d)_\CC$.
These results are useful for decomposing the Kuga Satake variety
into simple abelian subvarieties.

\newsubsection{K3 type Hodge structures with real multiplication.}
Assume now that $V$ is of K3 type with
$F=End_{Hod}(V)$ is a real field.
Then
$$
Hdg(V)(\CC)\cong SO(m,\CC)^n.
$$
One would again like to know 
$End(A)_\QQ={\rm End}_{Hdg(A)}(H^1(A,\QQ))$,
but this seems rather hard.
As a first step, we consider the decomposition of 
the Spin representation of the Lie algebra $so(d)$ 
upon restriction to $so(m)^n$, where $d=mn$.

Let $V_i$, $1\leq i\leq n$, be representations of $so(m)$, 
then we write $V_1\boxtimes V_2\boxtimes\ldots\boxtimes V_n$
for the representation of $so(m)^n$ where the $i$-th component
of $so(m)^n$ acts on $V_i$.

\newsubsection{The Spin representation.}
Let $S(m)$ be the Spin representation of the complex 
Lie algebra $so(m):=so(m)_\CC$ (cf.\ \cite{FH}, Chapter 20).
In case $m=2m'+1$ is odd, the Spin representation is an irreducible
representation of $so(m)_\CC$ of dimension $2^{m'}$.
In case $m=2m'$ is even, the Spin representation is
the direct sum of two irreducible components $S^\pm(2m')$:
$$
S(2m')=S^+(2m')\,\oplus\,S^-(2m'),\qquad \dim S(2m')=2^{m'},\quad
\dim S^\pm(2m')=2^{m'-1}.
$$
%We will omit the subscript $\CC$ from the Lie algebras in this section.
%For a non-degenerate quadratic form $Q$, 
%the $2^{n-1}$-dimensional complex even Clifford algebra $C^+(Q)_\CC$ is a %$so(Q)$-representation and, with $so(Q)_\CC\cong so(m)$,
%$C^+(Q)_\CC\cong S(m)^{2^{m'-1}}$
%if $m=2m'$ and $C^+(Q)_\CC\cong S(m)^{2^{m'}}$ if $m=2m'+1$.

\newsubsection{Lemma.}\label{spin restriction}
The restriction of the Spin representation 
$S(nm)$ of $so(nm)$ to $so(m)^n$ is given as follows:
$$
S(nm)_{|so(m)^n}=S(m)\boxtimes\ldots\boxtimes S(m)
\qquad\mbox{if}\quad m\equiv 0\;(2).
$$
In case $m$ is odd, write $n=2n'$ if $n$ is even and $n=2n'+1$ if
$n$ is odd. Then we have:
$$
S(nm)_{|so(m)^n}=
\bigl(S(m)\boxtimes\ldots\boxtimes S(m)\bigr)^{2^{n'}}
\qquad\mbox{if}\quad m\equiv 1\;(2).
$$ 

\ts
We use the conventions from \cite{FH}, Chapter 20.
In case $m$ is even, so is $nm$ and the Lie algebras $so(nm)$ and 
$so(m)^n$ both have rank $nm/2=n(m/2)$. Thus we can assume
that they have the same Cartan algebra $\goth\cong \CC^{nm/2}$
and the same dual $\goth^*$ generated by $L_1,\ldots,L_{nm/2}$
The weights of $S(nm)$ are then (\cite{FH}, Propostion 20.5)
the $(\pm L_1\pm L_2\ldots\pm L_{nm/2})/2$, each with multiplicity
one. Any such weight is the sum, in a unique way, of
the $n$ weights $(\pm L_{am/2+1}\pm \ldots\pm L_{(a+1)m/2})/2$ where
$a=0,1,\ldots,n-1$. As these are the weights of $S(m)\boxtimes\ldots\boxtimes S(m)$, with the same multiplicity, one,
we get the result.

In case $m$ is odd, the Lie algebra $so(m)$ has rank $(m-1)/2$.
If $n$ is odd, then one has:
$$
rk(so(nm))=(nm-1)/2=n(m-1)/2+(n-1)/2=rk(so(m)^n)+n'.
$$
We can now assume that $L_{n(m-1)/2+1},\ldots,L_{n(m-1)/2+n'}$
are zero on the Cartan algebra of $so(m)^n$.
The weights of $S(nm)$ are as in the previous case (\cite{FH}, 
proof of Proposition 20.20), but in the restriction to $so(m)^n$
there are $2^{n'}$ weights which map to the same sum of weights.
The case $n$ is even is similar.
\qed

\newsubsection{} In the remainder of this section we discuss two examples of Kuga Satake varieties of K3 type Hodge structures with real multiplication. In the first
example we have $n=[F:\QQ]=2,m=\dim_FV=3$
and we assume moreover that $C^+(\psi)\cong M_4(E)$, the 
algebra of $4\times 4$ matrices with coefficients in
an imaginary quadratic field $E$ (cf.\ \cite{L}, \cite{vG}, 9.2).
In the other case we take $n=[F:\QQ]=3,m=\dim_FV=3$.

\newsubsection{Proposition.}\label{KS n=2,m=3}
Let $(V,h,\psi)$ be a K3 type Hodge structure with $End_{Hod}(V)\cong F$, a real quadratic field and with 
$d=\dim V=6$. 
%Let $\psi$ be the quadratic form on $V$ 
%defined by $\psi$ and 
Assume that $C^+(\psi)\cong M_4(E)$ 
for an imaginary quadratic field $E$.

Then the Kuga Satake variety $A$ of $V$ is an abelian variety of 
dimension $16$,
$$
A\,\isog\, B^4,\qquad D:=End(B)_\QQ,\qquad NS(B)\cong\ZZ,
$$
where $B$ is a simple abelian fourfold and 
$D$ is a definite quaternion algebra over $\QQ$ which contains $E$.
There are three copies of the Hodge structure $V$ in $H^2(B,\QQ)$:
$$
H^2(B,\QQ)\cong V^3\oplus W,\qquad Hom_{Hod}(V,W)=0.
$$ 

\ts
As $C^+(\psi)\subset End(A)_\QQ$, we get
$$
End(A)_\QQ\cong M_4(E),\quad \mbox{hence}\quad
A\cong B^{4},\quad End(B)_\QQ=E,\quad
H^1(B,\CC)\cong S(6)
$$
where $B$ is an abelian variety of dimension $4$
and we use the isomorphism of $so(6)$-representations
$C^+(Q)_\CC\cong S(6)^{4}$ (\cite{vG}, Example 6.6).

As $V$ has real multiplication by the real quadratic field
$F$, the Lie algebra of $Hdg(V)_\CC$ is the subalgebra $so(3)^2$
of $so(6)$. 
Proposition \ref{spin restriction} shows that the restriction
of the Spin representation $S(6)$ of $so(6)$
to $so(3)^2$ is $(S(3)\boxtimes S(3))^2$.
The Spin representation
$S(3)$ of $so(3)\cong sl(2)$ is well-known to be the standard
two dimensional representation $V_1$ of $sl(2)$. 
We write $V_n$ for the irreducible representation of highest weight 
$n$ of $sl(2)$, note that $\dim V_n=n+1$. Thus we have an isomorphism
of $sl(2)^2$-representations:
$$
H^1(B,\CC)\cong \bigr( V_1\boxtimes V_1\bigr)^2.
$$

Next we compute $NS(B)_\CC=(\wedge^2 H^1(B,\CC))^{Hdg(B)}$,
using that $sl(2)^2$ is the Lie algebra of $Hod(B)_\CC$.
Note the following isomorphisms of $sl(2)\times sl(2)$-representations
$$
\begin{array}{rcl}
\wedge^2H^1(B,\CC) &\cong &\wedge^2((V_1\boxtimes V_1)^2)\\
&\cong&\left(\wedge^2(V_1\boxtimes V_1)\right)^2\,\oplus\,
(V_1\boxtimes V_1)\otimes(V_1\boxtimes V_1)\\
&\cong&(\wedge^2(V_1\boxtimes V_1))^3\,\oplus\, 
Sym^2(V_1\boxtimes V_1).
\end{array}
$$
It is not hard to check (using weights for example) that
$$
Sym^2(V_1\boxtimes V_1)\cong V_0\boxtimes V_0\,\oplus\, V_2\boxtimes V_2,
\qquad
\wedge^2(V_1\boxtimes V_1)\cong V_2\boxtimes V_0
\oplus V_0\boxtimes V_2,
$$
note that $V_0$ is the trivial representation of $sl(2)$.
In particular, there is a unique invariant in $H^2(B)$, 
so $NS(B)\cong \ZZ$ and hence $B$ is simple.

As $H^1(B,\CC)$ is the direct sum of two copies of an irreducible $Hdg(B)_\CC$-representation, it follows that 
$$
End(B)_\CC=End_{Hdg(B)}(H^1(B,\CC))\cong M_2(\CC),
$$
thus $D=End(B)_\QQ$ is a (non-commutative) division algebra of degree
four over $\QQ$, which contains the imaginary quadratic field $E$. 
To see that $D$ is definite, we must to show that $D_\RR$
is not isomorphic to $M_2(\RR)$. 
As elements of $D_\RR$ are endomorphisms
of $H^1(B,\RR)$ which commute with $Hdg(B)(\RR)$, it suffices
to show that $H^1(B,\RR)$ is an irreducible representation of $Hdg(B)(\RR)$. From \ref{real Hodge} we know that $Hdg(V)(\RR)\cong SO(2,1)
\times SO(3,\RR)$. 
The Spin group is then $SL(2,\RR)\times SU(2)$
and $H^1(B,\RR)$ is the $\boxtimes$-product of the standard
two dimensional representation of $SL(2,\RR)$ and the standard
two dimensional complex representation of $SU(2)$, 
which is an eight dimensional representation, irreducible over $\RR$. 
% quaternions are the commutant os su(2), and thus are ontained
%in the commutant of Sl(20\boxtimes su(2), but this commutant has 
% dim 4 over \RR and thus is the quaternions
% note su(2)\boxtimes su(2) is reducible!

The Hodge structure $V$ corresponds to the $Hdg(B)$-representation
with complexification:
$$
V_\CC=V_2\boxtimes V_0\oplus V_0\boxtimes V_2,
$$
thus $V^3$ is a Hodge substructure of $H^2(B,\QQ)$. 
\qed

\newsubsection{The case $m=n=[F:\QQ]=3$.}\label{fede}
For a general $h$ (so $End_{Hod}(V)=\QQ$) we have the following decomposition, up to isogeny, of the Kuga Satake variety $A$ of $V$:
$$
A\,\isog\, B^8,\qquad \dim B=16,\qquad End(B)_\QQ\cong D,\qquad 
H^1(B,\CC)\cong S(9)^2,
$$
where $D$ is a quaternion algebra
(cf.\ \cite{vG}, Proposition 7.7) and we used the isomorphism of
$so(9)$-representations $C^+(\psi)_\CC\cong S(9)^{16}$.
In case $D\cong M_2(\QQ)$, $B$ is isogenous to $B_1^2$ for 
an abelian eightfold $B_1$.

Assume now that $End_{Hod}(V)=F$, a totally real cubic extension of $\QQ$, so 
$Hdg(V)(\CC)\cong SO(3,\CC)^3$. Then the Lie algebra of the 
complex Hodge group of the Kuga Satake variety reduces from
$so(9)$ to $so(3)^3\cong sl(2)^3$ and 
the Spin representation $S(9)$ of $so(9)$ restricts
to $2$ copies of an $8$ dimensional irreducible representation 
(cf.\ Proposition \ref{spin restriction}, \cite{Ga}, Prop.\ 4.9, notation as in the proof of Proposition \ref{KS n=2,m=3})
$$
S(9)_{|sl(2)^3}\,\cong\,W^2,\qquad
W:=V_1\boxtimes V_1\boxtimes V_1.
$$
In particular, $End_{Hod}(H^1(B,\QQ))$ is a $\QQ$-algebra of rank
$4\cdot 4=16$. 

Assume that $End_{Hod}(H^1(B,\QQ))\cong M_4(\QQ)$.
Then we have an isogeny
$$
B\,\isog\, B_2^4,\qquad \dim B_2=4,\qquad End(B_2)_\QQ=\QQ,
$$
and $H^1(B_2,\CC)\cong W$. 
This case is discussed in \ref{fede mum}. The K3 type Hodge structure $V$ is not a Hodge substructure of $H^2(B_2,\QQ)$, 
but $H^2(B_2^2,\QQ)\cong V\oplus V'$ with $Hom_{Hod}(V,V')=0$
(\cite{Ga}, proof of Thm 3.4).

%\newpage

\section{The Kuga Satake variety and corestriction of algebras}
\label{corestriction}

\newsubsection{}
In this section we use the corestriction of algebras
to describe the first cohomology group of the
Kuga Satake variety as $Hdg(V)$-representation 
where $V$ is a K3 type Hodge structure with real multiplication. 
In contrast to the previous section, where we restricted the complex Spin representation to $Hdg(V)$, we now obtain a direct construction
over the rational numbers. This construction generalizes the
results of Galluzzi, cf.\ \ref{fede mum}, \ref{fede} and \cite{Ga}, which show that certain abelian varieties constructed by Mumford are Kuga Satake 
varieties.

\newsubsection{The corestriction.}\label{cores}
We use the notation from \ref{split F galois}:
$\tF$ is the Galois closure of a finite extension $F$ of $\QQ$, 
%as in \ref{split F galois}, 
$H=Gal(\tF/F)\subset G=Gal(\tF/\QQ)$ 
%are the Galois groups defined there 
and $n=[F:\QQ]=[G:H]$.
A coset $gH\in G/H$ gives a well-defined embedding
$F\hookrightarrow \tF$, $a\mapsto g(a)$ and thus defines
an $F$-algebra structure on $\tF$, this $F$-algebra
is denoted by $\tF_g$.

For an $F$-algebra $R$ and a coset $gH\in G/H$
the twisted algebra $R_g=R_{gH}$ is defined to be 
%the ring $R$ with $F$-algebra structure defined by 
(cf.\ \cite{M}, cf.\ \cite{S}, 8.8, \cite{R}, 4.4)
$$
R_g:=R\otimes_F\tF_g\qquad\mbox{so}\quad
ar\otimes 1=r\otimes g(a) \qquad (a\in F,\; r\in R).
$$
For $g\in G$ we have natural $\QQ$-linear maps:
$$
g:R_{g'}\longrightarrow R_{gg'},\qquad 
r\otimes a\longmapsto r\otimes g(a).
$$
To see that the map is well-defined, write
$a=g'(b)$, then 
$$
r\otimes a=br\otimes 1\longmapsto 
r\otimes g(a)=r\otimes (gg')(b)=br\otimes 1.
$$
Thus we get an action of $G$ on 
$$
Z_R:=
R_{g_1}\otimes_\tF\ldots\otimes_\tF R_{g_n},
\qquad\mbox{with}\quad
G/H=\{g_1H,\ldots,g_nH\}.
$$
The corestriction of the
$F$-algebra $R$ is the $\QQ$-algebra 
(\cite{M}, cf.\ \cite{S}, 8.9, \cite{R}, Thm 11, 5.5)
of $G$-invariants in $Z_R$:
$$
cores_{F/\QQ}(R)=Z_R^{G}.
$$
One has
$cores_{F/\QQ}(R)\otimes_\QQ \tF\cong Z_R$
and $\dim_\QQ cores_{F/\QQ}(R)=(\dim_F R)^n$.

Let $R^*$ be the multiplicative group of invertible elements of $R$.
Then there is a natural diagonal homomorphism:
$$
R^*\longrightarrow (cores_{F/\QQ}(R))^*,\qquad
u\longmapsto u\otimes u\ldots\otimes u,
$$
where $(cores_{F/\QQ}(R))^*$ is the multiplicative group of units
in the $\QQ$-algebra $cores_{F/\QQ}(R)$.

\newsubsection{Proposition.}\label{prop cores}
Let $(V,h,\psi=tr(\Phi))$ be a K3 type Hodge structure 
with $F=End_{Hod}(V)$ a totally real field.
Let $C^+(\psi)$ be the Kuga Satake Hodge structure associated to $V$
and let $C^+_F(\Phi)$ be the even Clifford algebra, over $F$,
of the $F$-bilinear form $\Phi:V\times V\rightarrow F$. 

Then $Hdg(C^+(\psi))$ is a subgroup of $(cores_{F/\QQ} C^+_F(\Phi))^*$
and there is an injective map of  
$Hdg(C^+(\psi))\,(\subset CSpin_F^+(\Phi))$ 
representations:
$$
cores_{F/\QQ} C^+_F(\Phi)\hookrightarrow C^+(\psi).
$$

\ts
We first extend scalars from $\QQ$ to $\tF$.
The Clifford algebra $C(\psi)$ of $\psi$ is the quotient
of the tensor algebra of $V$ by the two-sided ideal generated by the
$v\otimes v-\psi(v,v)$ for $v\in V$.
Extending scalars, we get an isomorphism 
$C(\psi)_\tF\cong C_\tF(\psi_\tF)$ where $\psi_\tF$
is the $\tF$-bilinear extension of $\psi$ to $V_\tF\times V_\tF$.

There is a direct sum decomposition of spaces 
with bilinear forms over $\tF$
(cf.\ the proof of Lemma \ref{dec psi}):
$$
(V_\tF,\psi_\tF)=\oplus_{g\in G/H}\,(V_g,\psi_g).
$$
This direct sum decomposition gives an isomorphism
(cf.\ \cite{S}, 9.2.5):
$$
C_\tF(\psi_\tF)=\,C_\tF(\psi_{g_1})\hat{\otimes}_\tF\ldots
\hat{\otimes}_\tF C_\tF(\psi_{g_n})
\qquad (G/H=\{H=g_1H,\ldots,g_nH\}),
$$
where $\hat{\otimes}$ is a graded tensor product.
It is easy to see that
$C^+_\tF(\psi_g)\cong C^+_\tF(\psi_e)_g$,
as $F$-algebras.
The weighted tensor product $\hat{\otimes}$ 
is the usual tensor product
on the `all even' part. So we have 
$$
C^+_\tF(\psi_{g_1})\otimes_\tF\ldots
\otimes_\tF C^+_\tF(\psi_{g_n})\cong
C^+_\tF(\psi_{e})\otimes_\tF\ldots
\otimes_\tF C^+_\tF(\psi_e)_{g_n}
\;\hookrightarrow\; C^+_\tF(\psi_\tF).
$$
As $V_e=V\otimes_F\tF_e$ and $\psi_e$ is the $\tF$-linear extension
of $\Phi$ to $V_e$, we get
$$
C^+_\tF(\psi_{e})\cong C^+_F(\Phi)\otimes_F\tF_e,\qquad
\mbox{hence}\quad
Z_{C^+_F(\Phi)}\hookrightarrow C^+_\tF(\psi_\tF)=C^+(\psi)_\tF.
$$
Taking $G$-invariants one finds that 
$cores_{F/\QQ}(C^+_F(\Phi))\hookrightarrow C^+(\psi)$.

The Hodge group of the weight one Hodge structure $C^+(\psi)$ 
is an algebraic subgroup of 
$CSpin(\psi)\subset C^+(\psi)^*$ and acts, by multiplication, on $C^+(\psi)$ and $C^+(\psi)_\tF$:
$u\cdot c:=uc$
for $u\in C^+(\psi)^*,c\in C^+(\psi)$ or $C^+(\psi)_\tF$.
Under the natural homomorphism $CSpin(\psi)\rightarrow GL(V)$,
$Hdg(C^+(\psi)$ maps onto $Hdg(V)=SO(V,\Phi)\subset SO(V)$.
In particular, $Hdg(C^+(\psi))\hookrightarrow CSpin_F(\Phi)$. 
The group $CSpin_F(\Phi)$ is a subgroup of $C^+_F(\Phi)^*$
which again acts by multiplication on $C^+_F(\Phi)$. 
Upon extending scalars
to $\tF$, it acts on $C^+_\tF(Q_{e})\cong C^+_F(\Phi)\otimes_F\tF_e$.
In particular, $Hdg(C^+(\psi))$ acts diagonally on
$Z_R$, with $R={C^+_F(\Phi)}$ and this gives the inclusion
$Hdg(C^+(\psi))\subset (cores_{F/\QQ} C^+_F(\Phi))^*$. 
This action is the restriction of the action of $Hdg(C^+(\psi))$ on $C^+(\psi)_\tF$. 
\qed

\newsubsection{Example.}  \label{fede mum}
Let $(V,h,\psi)$ be a K3 type Hodge structure with $End_{Hod}(V)=F$ a totally real field with $[F:\QQ]$=3 and assume that $\dim_FV=3$.
Then $\psi=tr(\Phi)$ and as $\Phi$ is defined on a three dimensional
$F$-vector space,  $C^+_F(\Phi)$ is a quaternion algebra $D$ over $F$ 
(cf.\ \cite{vG}, 7.5). 

As $\Phi_\sigma$ is indefinite  for one embedding and positive definite for the other two embeddings $\sigma:F\hookrightarrow \RR$ 
(cf.\ the proof of \ref{real Hodge}),
the algebra $D$ splits
for one embedding of $F\hookrightarrow \RR$
and is isomorphic to the quaternions for the other two embeddings,
so $D_\RR\cong M_2(\RR)\times \HH\times\HH$.

Conversely, a quaternion algebra defines a quadratic 
space $(D_0,N)$ over $F$, with $D_0$ the subspace of $D$ of elements with trace zero
and $N:D_0\rightarrow F$ the restriction of the norm on $D$ to $D_0$.
If $D_\RR$ is as above, then one can define K3 type Hodge structures on $V=D_0$ 
with endomorphism algebra $F$
(cf.\ the proof of Lemma \ref{no m=2}).

Assume that $cores_{F/\QQ}(D)\cong M_8(\QQ)$. 
Then Proposition \ref{prop cores} shows that
$Hdg(C^+(\psi))$ is a subgroup of $GL(8,\QQ)$. 
Thus the Kuga Satake variety of $V$ has a four dimensional abelian subvariety.
Actually, cf.\  \cite{M}, \cite{Ga},
$$
Hdg(A)=\ker(N:D^*\rightarrow F^*),\qquad A\,\isog\, B_2^{32},
$$ 
with $B_2$ a four dimensional abelian variety.

Mumford \cite{M} discovered these abelian varieties
and he showed that $End(B_2)_\QQ\cong \QQ$ but that nevertheless  
$Hdg(B_2)\neq Sp(8,\QQ)$. 
The relation with Kuga Satake varieties was established in \cite{Ga}.

\section{Predictions from the Hodge conjecture}

\newsubsection{The Hodge conjecture.}
The rational cohomology groups
$H^k(X,\QQ)$ of a smooth projective variety $X$ 
have a (polarized) rational Hodge structure of weight $k$.
The Hodge conjecture asserts that the space of
codimension $p$ Hodge classes
$$
B^p(X):=H^{2p}(X,\QQ)\cap H^{p,p}(X)
$$
is spanned by classes of algebraic cycles. The conjecture is
still very much open for $p\neq 0,1,\dim X-1,\dim X$.

\newsubsection{Hodge classes on a product.}
Let $X,Y$ be smooth projective varieties. The K\"unneth formula
and Poincar\`e duality $H^k(X,\QQ)\cong H^{2d_X-k}(X,\QQ)^{dual}$
imply that:
$$
H^k(X\times Y,\QQ)\,\cong\,
\oplus_{l+m=k} H^l(X,\QQ)\otimes H^m(Y,\QQ)\;\cong\,
\oplus_{l+m=k} Hom(H^{2d_X-l}(X,\QQ),H^m(Y,\QQ)).
$$
The summands $H^l(X,\QQ)\otimes H^m(Y,\QQ)$ are Hodge substructures
of $H^k(X\times Y,\QQ)$. The Hodge cycles in this summand are
exactly the homomorphisms of Hodge structures:
$$
B^p(X\times Y)\,\cap Hom(H^{2d_X-l}(X,\QQ),H^m(Y,\QQ))
\;=\;Hom_{Hod}(H^{2d_X-l}(X,\QQ),H^m(Y,\QQ)),
$$
where $2p=l+m$.

\newsubsection{Products of K3 surfaces.}
Let $X=S_1$, $Y=S_2$ be (algebraic) K3 surfaces. 
We consider the Hodge classes
in $H^4(S_1\times S_2)$. 
Note that $H^1(S)=H^3(S)=0$ for a K3 surface $S$. 
The summands $Hom_{Hod}(H^0(S_1),H^0(S_2))$, $Hom_{Hod}(H^4(S_1),H^4(S_2))$
are obviously spanned by the classes of $\{pt\}\times S_2$  
and $S_1\times\{pt\}$ respectively.

Recall that the Hodge structure on $H^2$ splits (cf.\ \ref{trans k3})
$$
H^2(S_i,\QQ)=NS(S_i)_\QQ\oplus T_{S_i,\QQ},
$$
and, as $NS(S_i)_\QQ^{2,0}=0$ and $T_{S_i,\QQ}$ is simple,
$$
Hom_{Hod}(NS(S_1)_\QQ,T_{S_2,\QQ})=0,\qquad
Hom_{Hod}(T_{S_1,\QQ},NS(S_2)_\QQ)=0.
$$
The vector space
$Hom_{Hod}(NS(S_1)_\QQ,NS(S_2)_\QQ)$
is spanned by classes of products of curves 
$C_1\times C_2\subset S_1\times S_2$.
Thus there remains the summand
$$
Hom_{Hod}(T_{S_1,\QQ},T_{S_2,\QQ}).
$$

\newsubsection{Hodge isometries.} 
Let $S$ be a K3 surface.
Then the Hodge structure $T_{S,\QQ}$ comes with a polarization
$\psi_S$, induced by the cup product on $H^2(S)$. 
A homomorphism of Hodge structures 
$$
f\in Hom_{Hod}(T_{S_1,\QQ}, T_{S_2,\QQ}),\qquad
\mbox{such that}\quad
\psi_{S_2}(f(v),f(w))=\psi_{S_1}(v,w)
$$
for all $v,w\in T_{S_1,\QQ}$ is called a Hodge isometry.

Mukai has announced that 
if $f\in Hom_{Hod}(T_{S_1,\QQ},T_{S_2,\QQ})$ is a Hodge isometry, then
$f$ is the class of an algebraic cycle on $S_1\times S_2$
(\cite{Mukai}, Theorem 2).
Under certain conditions
on the dimension of $T_{S_1,\QQ}$ proofs were given earlier by Mukai and Nikulin (cf.\ \cite{Mukai}, section 4).
This solves the Hodge conjecture for Hodge isometries, 
but below we recall that there is still a lot to do.

\newsubsection{Complex multiplication.}
In case $S_1=S_2$ and $End_{Hod}(T_{S_1,\QQ})$ is a CM field,
Ramon Mari \cite{pepe} showed that $End_{Hod}(T_{S_1,\QQ})$ is 
spanned by Hodge isometries. Thus, by Mukai's results, 
any $f\in End_{Hod}(T_{S_1,\QQ})$ is the class of an algebraic cycle.

\newsubsection{Real multiplication.}
In case $S_1=S_2$, the rational multiples of the identity can
be obtained from the projection to $End_{Hod}(T_{S_1,\QQ})$ of 
rational multiples of the class of the diagonal in $S_1\times S_1$.

However, if $End_{Hod}(T_{S_1,\QQ})$ is a totally real field, distinct
from $\QQ$, I do not know of any example where a non-trivial
endomorphism is represented by an algebraic cycle.

\newsubsection{Scaling the polarization.}
Let $(T_{S,\QQ},\psi_S)$ be the polarized Hodge structure
defined by a K3 surface $S$. For any positive integer $n$,
there is the polarized Hodge structure of K3 type 
$(T_{S,\QQ},n\psi_S)$. 

In general, these Hodge structures are not
Hodge isometric. For example, if $d=\dim T_{S,\QQ}$ is odd, and $n$
is not a square, then as $\det(n\psi_S)= n^d\det(\psi_S)$, we get
an obstruction to the existence of isometries (cf.\ \ref{det and dis}).

This construction is a special case of the real twist of \ref{twist k3},
$n\psi=\psi_n$ for $a=n\in\ZZ$. Thus for $n\leq 11$
there exists a K3 surface $S_n$ such that
its transcendental lattice $(T_{S_n},\psi_{S_n)}$ is isometric to
$(T_S,n\psi_S)$. In particular, the identity on $T_S$ gives a non-trivial element in $Hom_{Hod}(T_{S,\QQ},T_{S_n,\QQ})$.

In case $n=2$, some of these homomorphisms of Hodge structures
can be shown to be the classes of algebraic cycles by using 
Nikulin involutions of K3 surfaces, cf.\ \cite{GL}, \cite{vGS}. 
In general, it seems to be an interesting open problem to find
such algebraic cycles.

\newsubsection{Twisting the polarization.}
Let $S$ be a K3 surface with 
$End_{Hod}(T_{S,\QQ})=F$ a totally real field.
Any totally positive $a\in F$ defines a polarization 
$(\psi_S)_a$ on $T_{S,\QQ}$ (cf.\ \ref{mult of pols}). 
In case $F=\QQ$ we recover the scaling operation described above.
In general, the polarized Hodge structures
$(T_{S,\QQ},\psi_S)$ and $(T_{S,\QQ},(\psi_S)_a)$ are
not isometric. 

If $d=\dim T_{S,\QQ}\leq 11$, there exists a K3 surface
$S_a$ such that $(T_{S_a},\psi_{S_a})$  is Hodge isometric
to $(T_S,(\psi_S)_a)$, cf.\ \ref{twist k3}. 
The identity map $T_{S,\QQ}\rightarrow T_{S,\QQ}=T_{S_a,\QQ}$
is a homomorphism of Hodge structures and, again, 
it seems to be an interesting open problem to show that is induced
by an algebraic cycle.

\newsubsection{Remark.}
Let $S_1,S_2$ be K3 surfaces. 
Let $Z$ be a smooth surface, with maps
$$
S_1\;\stackrel{\pi_1}{\longleftarrow}\;Z\;
\stackrel{\pi_1}{\longrightarrow}\;S_2\qquad
\mbox{let}\quad
f=\pi_{2*}\pi_1^* \in Hom_{Hod}(T_{S_1,\QQ},T_{S_2,\QQ})
$$
be the homomorphism of Hodge structures defined by
the class of $(\pi_1\times\pi_2)(Z)\subset S_1\times S_2$.
Assume that $f\neq 0$.
Let $V_i:=T_{S_i,\QQ}$, as the $V_i$ are simple and
$f\neq 0$ we get an isomorphism of rational Hodge structures 
$f:V_1\stackrel{\cong}{\longrightarrow} V_2$.

The map $f$ is not necessarily an isometry $(V_1,\psi_1)\rightarrow (V_2,\psi_2)$, where $\psi_i:=\psi_{S_i}$, the polarization induced
by cupproduct on the surface $S_i$. 
In fact, $f=\pi_{2*}\pi_1^*$ and as
$\pi_1^*:V_1\hookrightarrow H^2(Z,\QQ)$ 
is compatible with the cupproduct:
$$
\psi_Z(\pi_1^*x_1,\pi_1^*y_1)
:=\pi_1^*x_1\cup \pi_1^*y_1
=\pi_1^*(x_1\cup y_1)=d_1(x_1\cup y_1)=d_1\psi_1(x_1,y_1),
$$
where $d_1$ is the degree of $\pi_1$ (that is, the cardinality of a general fiber). So we have Hodge isometries 
$$
(V_1,\psi_1)\,\longrightarrow\, (V_1,d_1\psi_1)\;\hookrightarrow\; (H^2(Z,\QQ),\psi_Z).
$$
The map $\pi_{2*}$ is not compatible with cupproducts, 
but the projection formula 
%$x_2\cup \pi_{2*}y_z=\pi_2^*x_2\cup y_z$, so
gives:
$$
\psi_{2}(x_2,\pi_{2*}y_z)=\psi_Z(\pi_2^*x_2,y_z)\qquad (x_2\in H^2(S_2,\QQ),\;y_z\in H^2(Z,\QQ)).
$$
%(we used that $\pi_{2*}:H^4(Z,\QQ)
%\stackrel{\cong}{\longrightarrow} H^4(S_2,\QQ)$

Assume now that $H^2(Z,\QQ)= \pi_1^*(V_1)\oplus W$ with $Hom_{Hod}(V_1,W)=0$,
so there is a unique copy of the Hodge structure $V_1\cong V_2$ in $H^2(Z,\QQ)$.
The composition $\pi_{2*}\pi_2^*$ is multiplication by $d_2$,
the degree of the map $\pi_2$, on $H^2(S_2,\QQ)$.
Thus the map $\pi_2^*:V_2\stackrel{\cong}{\longrightarrow} \pi_1^*(V_1)$ is an isomorphism and
%$\pi_2^*:V_2\stackrel{\cong}{\longrightarrow} \pi_1^*(V_1)$
%with inverse $d_2^{-1}\pi_{2*}$.
%Therefore, as before, $\psi_Z(\pi_2^*x_2,\pi_2^*y_2)=d_2\psi_2(x_2,y_2)$.
given $x_z,y_z\in f(V_1)\subset H^2(Z,\QQ)$, there are $x_2,y_2\in V_2$ with $x_z=\pi_2^*x_2,y_z=\pi_2^*y_2$. Therefore:
$$
\psi_2(\pi_{2*}x_z,\pi_{2*}y_z)=
\psi_2(\pi_{2*}\pi_2^*x_2,\pi_{2*}\pi_2^*y_2)=d_2^2\psi_2(x_2,y_2)
=d_2\psi_Z(\pi_2^*x_2,\pi_2^*y_2)
=d_2\psi_Z(x_z,y_z).
$$
%from which we get:
%$$
%\psi_2(\pi_{2*}x_z,\pi_{2*}y_z)=d_2\psi_Z(\pi_2^*x_2,\pi_2^*y_2).
%$$
In particular, the isomorphism $f:V_1\rightarrow V_2$ induces a scaling
on the polarizations:
$$
\psi_2(fx_1,fx_2)=d_1d_2\psi_1(x_1,x_2).
$$

In general, given a codimension two cycle $\sum n_iZ_i$ on $S_1\times S_2$, 
after replacing each surface $Z_i$ by its desingularization
(which maps to $S_1\times S_2$ with image $Z_i$), the homomorphism
of Hodge structures induced by $\sum n_iZ_i$ is a linear combination
of maps as above. Thus to get an `interesting' map 
$f:V_1\rightarrow V_2$, 
induced by an algebraic cycle,
one needs surfaces $Z_i$ whose $H^2$ 
contains more than one copy of $V$.

\newsubsection{The Kuga Satake Hodge conjecture.}
The Kuga Satake variety $A$ of a K3 type Hodge structure $(V,h,\psi)$
has the property that there is a homomorphism of Hodge structures
$V\hookrightarrow H^2(A^2,\QQ)$ (cf.\ \cite{vG}, 6.3.3).
In particular, if $V\cong T_{S,\QQ}\subset H^2(S,\QQ)$ 
for a K3 surface $S$, then
the Hodge conjecture predicts the existence of a cycle $Z$ on
$S\times A\times A$ which induces an isomorphism from the copy of $V$
in $H^2(S)$ to the one in $H^2(A^2)$ .

If there is such a cycle $Z$ and
if $End_{Hod}(T_{S,\QQ})$ is also generated by algebraic cycles
then any homomorphism of Hodge structures
$T_{S,\QQ}\rightarrow  V\subset H^2(A^2,\QQ)$ 
is represented by an algebraic cycle on $S\times A^2$.

\

%Version: \today $\;$ k3rm.tex
%\bibliographystyle{amsplain}
%\bibliography{biblio}
\end{document}